\documentstyle[12pt,twoside]{article}

\newtheorem{theorem}{Theorem}

\begin{document}

\overfullrule=0pt
\baselineskip=24pt
\font\tfont= cmbx10 scaled \magstep3
\font\sfont= cmbx10 scaled \magstep2
\font\afont= cmcsc10 scaled \magstep2
\title{\tfont Intersection Forms and the Adjunction Formula for Four-manifolds
via  CR Geometry }
\bigskip
\author{ Mikhail Chkhenkeli\\ Thomas Garrity \\
\\ Department of Mathematics\\ Williams College\\Williamstown, MA  
01267\\ 
email: Mikhail.Chkhenkeli@williams.edu\\
 tgarrity@williams.edu}
\date{}
\maketitle
\begin{abstract}

This is primarily an expository note showing that
 earlier work of Lai\cite{Lai} on CR geometry provides a clean 
 interpretation, in terms of a Gauss map,
  for an adjunction formula for embedded surfaces in 
an almost complex four manifold.
 We will see that if $F$
 is a surface with genus $g$ in an almost complex four-manifold $M$, then 
	  $$  	2 - 2 g + F \cdot F  - i^{*} c_{1}(M) - 2 F\cdot C = 0,$$
where $C$ is a two-cycle on $M$ pulled back from the  cycle of 
two planes with complex structure in a Grassmannian $Gr (2, C^{N})$ via
 a Gauss map and where $i^{*} c_{1}(M)$ is the restriction of the first
  Chern class of $M$ to $F$.    The key new term of interest is $F 
  \cdot C$, 
  which will capture the points of $F$ whose tangent planes inherit a
   complex structure from the almost complex structure of the ambient 
   manifold $M$.  These complex 
   jump points then determine the genus
    of smooth representatives of a homology class in  $H_{2}(M, Z)$.  
    Further, via polarization, we can use this formula to determine 
    the intersection form on $M$ from knowing the nature of the 
    complex jump points of $M$'s surfaces.

\end{abstract}

\section{The Adjunction Formula}

This goal of this note is primarily to point out how a paper of Lai 
from the early 70s can be used to interpret the intersection form of 
an almost complex four-manifold.  This was earlier pointed 
out by Eliashberg and Harlamov\cite{EH}, \cite{E}.

Let $M$ be an almost complex four-manifold.  Then  at each point
$ p$ of $M$ there is an automorphism 
$J: T_{p}M \rightarrow T_{p}M$  with   $J^{2} = -I$   and such that
$  J$  varies smoothly on $M$.  Thus at each point $p$
 we can identify the tangent plane  $T_{p}M$ with the complex
 two-plane $C^{2}$.  Let $F$ be an embedded surface in $M$.
   At most points $p$ on $F$, we expect 
$$		                       T_{p}F \bigcap JT_{p}F = 0,$$
meaning that  $T_{p}F$  will not inherit a complex structure from 
$T_{p}M$ . But there will be some points at which 
		                          $$T_{p}F  = JT_{p}F,$$
namely those points whose real tangent plane inherit  the 
 structure of a complex line from  $T_{p}M$.
   We call these points {\it complex jump points} (In \cite{Lai}, Lai
    used the term RC-singular point and in \cite{Wells}, Wells used the term 
    nongeneric point). The complex jump points will provide a 
    sharp link between  the Euler characteristic of $F$, 
    the self-intersection number of $F$
     and the pullback of the first Chern class of $M$ to $F$.
       Namely, we will show the adjunction formula
$$	          2 - 2 g + F \cdot F  - i^{*} c_{1}(M) - 2 F \cdot C = 0,$$

where $i: F \rightarrow  M$  is the embedding map and $F \cdot C$
 represents, as we will see, the complex jump points.

          We first saw such a formula in a lecture by 
          Kirby and realized that there was a natural
           proof via the Gauss map.  This led us to
            the fact that the Kirby formula actually 
            was a special case of an extension of an
             old formula of Lai, who was working in CR geometry. 
              Kirby then pointed out to us the earlier work of 
               Eliashberg and Harlamov, for which we thank him.

The key will be the Gauss map, which is the natural map 

$$		                    \sigma: F \rightarrow GR(2, C^{N}),$$
given via the embedding  $i: T_{p}F \rightarrow T_{p}M$
and then choosing enough sections of $TM$
 to have a mapping into a complex affine space.  Set
$$		C  := (\Lambda \in  Gr(2, C^{N}): J \Lambda   = \Lambda ), $$
where $J$ is the automorphism associated with the complex structure on
$ C^{N}$.  The homology of  $C$
 can be explicitly computed in terms of the special Schubert cycles of  
 $Gr(2, C^{N})$ . 
The complex jump points of the surface $F$ are precisely the pullback of 
the points  $\sigma(F) \cdot C$.
	Lai showed in \cite{Lai}:
\begin{theorem}[Lai]
 Let $F$ be a compact real k-dimensional manifold and $M$ a real
  2n-dimensional almost complex manifold. Let $i: F \rightarrow  M$
   be an immersion. Assume $2n-2=k$. Then
$$\Omega(F) + \sum_{r=0}^{n-1}\bar{\Omega}(F)^{n-r-1}\cup i^{*}c_{r}(M) =
2\sigma ^{*}(\sigma (F) \cdot C).$$
\end{theorem}
Here $\Omega(F)$ is the Euler class of $F$,  $\bar{\Omega}(F)$
  is the Euler class of the normal bundle of $F$ in $M$ and    $\sigma$
    and  $C$
        are the higher dimensional analogues of our earlier definitions.
          For real surfaces F in complex surfaces (and 
          thus real four manifolds) M, Lai's formula becomes
$$\Omega(F) + \bar{\Omega}(F) -  i^{*}c_{1}(M) =
2\sigma^{*}(\sigma (F) \cdot C).$$
But in this case, the Euler characteristic of the normal bundle 
can be identified with the self-intersection number  $F \cdot F$
.  Using that the Euler characteristic of $F$ can be identified to  
$2 - 2 g$, we have at least formally our desired adjunction formula.  To
 actually prove the formula, all we need do is to examine Lai's proof,
  in which it is assumed that $ M $ is a complex manifold, and then 
  simply to observe that all he needed to use  was that $M$ has the
   structure of an almost complex manifold.  Thus 
   the real purpose of this paper is to point out that Lai's formula, which
    he no doubt developed without thinking at all about
     four-manifolds but instead 
     about low codimensional CR structures, can be easily applied
      to the geometry of four-manifolds.  

Note that we have altered the notation of Lai. His $M$ is our $F$,
 his $N$ is our $M$ and his $DK$ is our
$\sigma(F) \cdot C$.

 \section{Intersection Forms}
 
 The key  to the topology of a four-manifold lies in understanding the 
 intersection form on $H_{2}(M,{\bf Z})$.  The Kirby-Lai formula 
 gives us
 $$  	 F \cdot F  = 2g -2 +  i^{*} c_{1}(M) + 2 F\cdot C,$$
and thus expressing the self-intersection number of a surface $F$ in 
terms of its genus, the smooth structure of $M$ and the algebraic 
number of complex jump points of $F$.  Once we have the 
self-intersection number, by polarization we will be able to recover 
the intersection form, in terms of the geometry of the manifold $X$.

Let $F$ and $G$ be two elements of $H_{2}(M,{\bf Z})$.  With only a 
slight abuse of notation, we can let $F$ and $G$ denote smooth 
representatives of $F$ and $G$ and let $F+G$ denote a smooth 
representative of the homology class of $F+G$.
Then we have
$$F\cdot G = \frac{1}{2}(F+G)\cdot (F+G) - \frac{1}{2}F\cdot F 
-\frac{1}{2} G\cdot G$$
$$= \frac{1}{2}(2g_{(F+G)} -2 +  i_{(F+G)}^{*} c_{1}(M) +  2(F+G)\cdot C)
- \frac{1}{2}(2g_{(F)} -2 +  i_{F}^{*} c_{1}(M) +  2F\cdot C)$$
$$-\frac{1}{2}(2g_{(G)} -2 +  i_{G}^{*} c_{1}(M) +  2G\cdot C).$$
Thus the topology of the intersection form is indeed captured by the 
geometry of the genus of the three surfaces $F$, $G$ and $F+G$,
 the first Chern class of $M$ and the number of 
complex jump points of $F$, $G$ and $F+G$.

 Of course, in practice it 
is unlikely that one would know this information without already 
knowing the intersection form.  Also, note that we only need the 
geometric information about smooth representatives for a basis for 
$H_{2}(M,{\bf Z})$ and for smooth representatives for pairwise sums 
of basis elements.

\section{Bounds on the number of complex jump points on
 characteristic 2-spheres}

In this section we obtain bounds on the number
 of complex jump points on $F$ when $F$ is a 2-sphere
smoothly representing a characteristic class of
 a closed, oriented, simply connected smooth four-manifold $M$.
  By definition, a homology class in $H_{2}(M)$
   is called characteristic if it is dual to $w_{2}(M) mod2$. In 
   \cite{Chkhenkeli} bounds on the self-intersection number of a characteristic 
2-sphere are given. By using these bounds and the Kirby - Lai adjunction formula, we obtain 
the following can be shown.
\begin{theorem} Let  $i: F \rightarrow M$
 be an embedding of a 2-sphere $F$ representing a characteristic class 
of a closed, oriented, simply connected smooth four-manifold $M$.
 Let $n$ be the number of complex jump points on $F$. Then   
$$( b^{+} - 9b^{-} + 10 - i^{*}c_{1}(M) )/2 \leq  n \leq 
( b^{+} - b^{-}/9 - 10/9 - i^{*}c_{1}(M) )/2 ,$$
 if $ F\cdot F \leq -1$;
$$( b^{+}/9 - b^{-} + 26/9 - i^{*}c_{1}¥(M) )/2 \leq n \leq
 ( 9b^{+} - b^{-} - 6 - i^{*}c_{1}(M) )/2 ,$$
  if $F\cdot F \geq 1$.
  \end{theorem}
Here $b^{+}$ and $b^{-}$ are the numbers of
 positive and negative eigenvalues of the intersection form of $M$,
respectively; $F\cdot F$ is the self-intersection number of $F$
\ and $i^{*}c_{1}(M)$ is the pullback of the first Chern
class of $M$ to $F$.


\begin{thebibliography}{99}





\bibitem{Chkhenkeli}  M. Chkhenkeli, Characteristic 2-Spheres in 4-Manifolds, 
in preparation,
 Williams College, 1997.
 
 \bibitem{E}  Y. Eliashberg, Filling by holomorphic discs and its applications,
 {\bf Geometry of low-dimensional manifolds}, 2 (Durham, 1989), London 
 Math. Soc. Lecture Note Ser., 151, Cambridge Univ. Press, Cambridge, 
 1990, pp. 45-67.
 
 \bibitem{EH}  V. M. Harlamov and Y. Eliashberg, On the number of 
 complex points in a complex surface, {\it Proc. of Leningrad Int. 
 Topology Conference}, 1982, 143-148.
 
 

 
\bibitem{Kirby} R. Kirby, An Adjunction Formula for Smooth Surfaces in 4-Manifolds,
 in preparation.
 
 \bibitem{Lai}  H. F. Lai, Characteristic Classes of Real Manifolds Immersed in Complex 
Manifolds, {\it Transactions of the American Mathematical Society}, Vol. 172 (1972), 
1-33.


\bibitem{Wells} R. O. Wells, Jr., Holomorphic Hulls and Holomorphic Convexity, 
{\it Complex 
Analysis}
(Proc. Conf. Rice Univ., Houston, Tex., 1967), Rice Univ. Studies 54 (1968),
 no. 4, 75-84.











			








	
				

	





\end{thebibliography}
\end{document}